\documentclass[12pt]
{amsart}

\usepackage{amssymb}
\usepackage{amsmath}
\usepackage{array}

\newcommand{\Lg}{\mbox{$\mathfrak g$}}
\newcommand{\Ll}{\mbox{$\mathfrak l$}}

\newcommand{\La}{\mbox{$\mathfrak a$}}
\newcommand{\Lm}{\mbox{$\mathfrak m$}}

\newcommand{\Lt}{\mbox{$\mathfrak t$}}

\newcommand{\B}[1]{\mbox{$\mathbf{B}_{#1}$}}

\newcommand{\D}[1]{\mbox{$\mathbf{D}_{#1}$}}
\newcommand{\G}{\mbox{$\mathbf{G}_2$}}
\newcommand{\E}[1]{\mbox{$\mathbf{E}_{#1}$}}
\newcommand{\F}{\mbox{$\mathbf{F}_4$}}

\newcommand{\Pf}{{\em Proof}. }
\newcommand{\EPf}{\hfill$\square$}

\newcommand{\R}{\mbox{$\mathbf R$}}
\newcommand{\CC}{\mbox{$\mathbf C$}}

\newcommand{\cohom}[2]{\mbox{$\mathrm{cohom}(#1,#2)$}}
\newcommand{\homrk}[2]{\mbox{$\mathrm{homrk}(#1,#2)$}}
\newcommand{\rk}{\mbox{$\mathrm{rk}\;$}}

\newcommand{\SU}[1]{\mbox{$\mathbf{SU}(#1)$}}
\newcommand{\U}[1]{\mbox{$\mathbf{U}(#1)$}}
\newcommand{\SP}[1]{\mbox{$\mathbf{Sp}(#1)$}}
\newcommand{\SO}[1]{\mbox{$\mathbf{SO}(#1)$}}
\newcommand{\OG}[1]{\mbox{$\mathbf{O}(#1)$}}
\newcommand{\Spin}[1]{\mbox{$\mathbf{Spin}(#1)$}}

\newtheorem{thm}{Theorem}

\newtheorem{rmk}[thm]{Remark}
\newtheorem{cor}[thm]{Corollary}
\newtheorem{prop}[thm]{Proposition}
\newtheorem{lem}[thm]{Lemma}

\title{Homogeneity rank of real representations\\ of compact Lie groups}
\author{Claudio Gorodski and Fabio Podest\`a}
\address{Instituto de Matem\'atica e Estat\'\i stica\\Universidade
de S\~ao Paulo\\Rua do Mat\~ao, 1010\\
S\~ao Paulo, SP 05508-900\\
Brasil}
\email{gorodski@ime.usp.br}
\address{Dipartimento di Matematica e Appl.\ per l'Architettura\\
Piazza Ghiberti 27\\50100 Florence\\Italy}
\email{podesta@math.unifi.it}
\date{\today}

\begin{document}

\begin{abstract}
The main result of this paper is the classification of
the real irreducible representations of compact Lie groups
with vanishing homogeneity rank.
\end{abstract}

\maketitle

\section{Introduction}

Let a compact Lie group $G$ act smoothly on a
smooth manifold $M$. The codimension of the principal
orbits in $M$ is called the
cohomogeneity $\cohom GM$ of the action.
P\"uttmann, starting from an inequality for the dimension of the fixed point set
of a maximal torus in $G$ due to Bredon (\cite{B},~p. 194), introduced in~\cite{Pu} the
\emph{homogeneity rank} of $(G,M)$ as the
integer
\begin{eqnarray*}
\homrk GM & = & \rk G-\rk G_\mathrm{princ} - \cohom GM \\
          & = & \rk G - \rk G_\mathrm{princ} + (\dim G -\dim G_\mathrm{princ}) - \dim M,
\end{eqnarray*}
where $G_\mathrm{princ}$ is a principal isotropy subgroup of the
action and, for a compact Lie group $K$,
$\rk K$ denotes its rank, namely the dimension
of a maximal torus. We will see in the next section
that orbit-equivalent actions have the same homogeneity rank.

This invariant, although not with this name, had already been
considered by Huckleberry and Wurzbacher who proved that a Hamiltonian
action of a compact Lie group on a symplectic manifold has vanishing
homogeneity rank if and only if the
principal orbits are coisotropic with respect to the invariant symplectic
form (see~\cite{HW}, Theorem~3, p.~267 for this result and
other characterizations of this property).
If $\rho:G\to \U{V}$ is a complex representation where
$V$ is a complex vector space endowed with an invariant symplectic
structure, then the $G$-action is automatically Hamiltonian, and
it has vanishing homogeneity rank if and only if every principal orbit is coisotropic; this condition can
be proved to be equivalent to the fact that a Borel subgroup of the
complexified group $G^c$ has an
open orbit in $V$, and also to the fact that the naturally induced
representation of $G$ on the ring of
regular functions $\CC[V]$ splits into the sum of mutually
inequivalent irreducible representations
(see e.g.~\cite{Kr}, p.~199).
Complex representations with these equivalent properties
are called \emph{coisotropic} or
\emph{multiplicity-free}; Kac~\cite{Kac} classified the irreducible
multiplicity-free
representations and, later, Benson and Ratcliff~\cite{BR} and,
independently, Leahy~\cite{L} classified the reducible ones.

In this paper we consider the case of an irreducible
representation $\rho:G\to\OG V$ of a compact Lie group $G$ on a
real vector space $V$ with vanishing homogeneity rank.
Since representations
admitting an invariant complex structure have null homogeneity
rank if and only if they are multiplicity-free, we will deal only with
irreducible representations of real type, also called
absolutely irreducible, namely those which admit no invariant
complex structure. Our main result is the following theorem.

\begin{thm}\label{thm:main}
An absolutely irreducible representation $\rho$ of a
compact connected Lie group $G$ has vanishing homogeneity rank if and only
if it is either orbit-equivalent to the isotropy representation of a
non-Hermitian symmetric space of inner type or it is one of the
following representations:
\[\begin{array}{|c|c|c|c|}
\hline
 G & \rho & d & c \\
\hline
\SP1\times\SP{n},\ n\geq2 & {\operatorname{(standard)}}^3
\otimes_{\mathbf{H}} {\operatorname{(standard)}} & 8n & 3 \\
\SO4\times \Spin7 & {\operatorname{(standard)}}\otimes_{\mathbf{R}} {\operatorname{(spin)}} & 32 & 5 \\
\SP1\times \Spin{11} & {\operatorname{(standard)}}\otimes_{\mathbf{H}}{\operatorname{(spin)}} & 64 & 6\\
\hline
\end{array}\]
where $d$ denotes the dimension of the representation space and $c$
denotes its cohomogeneity.
\end{thm}

\section{Preliminaries}

Let $(G,V)$
be an absolutely irreducible representation of a compact Lie
group $G$ on a real vector space $V$.
It is shown in Corollary~1.2 in~\cite{Pu} that the
homogeneity rank of a linear representation is non positive. In
this regard, the representations with vanishing homogeneity rank
are precisely those with maximal homogeneity rank. The following
monotonicity property that is stated on p.~375 in~\cite{Pu} and is
valid for smooth actions on smooth manifolds will be the basis of
the method of our classification. Since there is no proof
in~\cite{Pu}, we include one for the sake of completeness.

\begin{prop}
Let $(G,M)$ be a smooth action.
If $G'$ be a closed subgroup of $G$, then $\homrk {G'}M\leq
\homrk GM$.
\end{prop}

\Pf We first prove the statement in the case in which $M$
is $G$-homogeneous, i.~e.~we prove that given a homogeneous
space $M=G/H$, where $G$ is a compact Lie group and $H$
is a closed subgroup, for every closed subgroup $G'$ of $G$ we have
\[ \homrk{G'}{G/H}\leq\rk G-\rk H. \]
We prove this by induction on the dimension of the manifold, the
initial case $\dim M=1$ being clear. Fix the point $o=[H]\in G/H$,
a maximal torus $T_H$ of $H$, and a maximal torus $T$ of $G$ containing
$T_H$. Since conjugation of $G'$ by elements of $G$ does
not affect the homogeneity rank, we can assume that a maximal torus
$T'$ of $G'$ sits inside $T$. Then we have
\[ \rk{G'}-\rk{G'_o}\leq\dim T'-\dim(T'\cap G'_o)=\dim T'\cdot o\leq\dim T\cdot o, \]
where $G'_o$ denotes the isotropy subgroup of $G'$ at $o$.
Therefore
\[ \rk{G'}-\rk{G'_o}\leq\rk G-\rk H. \]
We now consider the slice representation of $G'_o$ on the normal
space $W$ to the orbit $G'\cdot o$; we can assume that the
dimension $k$ of $W$ is at least $2$, since otherwise $G'_o$
contains a principal isotropy subgroup of $(G,G/H)$ as a subgroup
of finite index and the claim follows immediately. Denote by $S$
the unit sphere in $W$ with respect to a $G'_o$-invariant inner
product in $W$ and apply the induction hypothesis. Since $G'_o$ is
a closed subgroup of $\SO k$, we have
\begin{eqnarray*}
\homrk{G'_o}S & = & \rk G'_o-\rk G'_{\mathrm{princ}}-\cohom{G'}M+1\\
              & \leq & \homrk{\SO k}S \\
              & = & \frac{1+(-1)^k}2,
\end{eqnarray*}
where $G'_{\mathrm{princ}}$ denotes a principal isotropy
subgroup of $G'$ on $M$. It then follows that
\[ \homrk{G'}M\leq\rk{G'}-\rk{G'_o}+\frac{-1+(-1)^k}2\leq\rk G-\rk H, \]
and our claim is proved.

In the general case, we fix a $G$-regular point $p\in M$ and observe that
a point $q\in G\cdot p$ is principal for the $G'$-action on $G\cdot p$
if and only if it is principal for the $G'$-action in $M$; this means
that $\cohom{G'}M=\cohom GM+\cohom{G'}{G\cdot p}$. We know from
the previous case that
\[ \homrk{G'}{G\cdot p}\leq\homrk G{G\cdot p}, \]
and now our claim follows by subtracting $\cohom GM$ from both
members of the above inequality. \EPf

\begin{cor}\label{cor:s-repr}
Let $(G,V)$ be a representation of a compact Lie group $G$ on a
real vector space $V$. If $(G,V)$ is not of vanishing homogeneity rank,
then the action of a closed subgroup of $G$ on $V$
is never of vanishing homogeneity rank.
\end{cor}

The preceding corollary indicates a strategy to classify
representations with vanishing homogeneity rank. First we observe
that the standard representation of $\SO n$ on $\R^n$ is of
vanishing homogeneity rank if and only if $n$ is even. Then we need
to decide which of the maximal subgroups of $\SO n$, where $n$ is
even, act absolutely irreducibly on $\R^n$ with vanishing homogeneity rank.
For each example that we encounter, we examine which of its maximal
subgroups still act absolutely irreducibly
on $\R^n$ with vanishing homogeneity rank, and
so on. The process will eventually yield all the closed subgroups
of $\SO n$ that act absolutely
irreducibly on $\R^n$ with vanishing homogeneity rank.
The effectiveness of this strategy is elucidated by the following well
known result of Dynkin~\cite{D}.
\begin{thm}[Dynkin]\label{thm:Dynkin}
\begin{enumerate}
\item Let $G$ be a maximal connected subgroup of $\SO n$. Then $G$ is
conjugate in $\OG n$ to one of the following:
\begin{enumerate}
\item $\SO k\times\SO{n-k}$, where $1\leq k\leq n-1$;
\item $\rho(\SO p\times\SO q)$, where $pq=n$ and $3\leq p\leq q$,
and $\rho$ is the real tensor product of the vector representations;
\item $\U k$, where $2k=n$;
\item $\rho(\SP p\times\SP q)$, where $4pq=n\neq4$,
and $\rho$ is the quaternionic tensor product of the vector representations;
\item $\rho(G_1)$, where $G_1$ is simple and $\rho$ is a real form of
a complex irreducible representation of degree $n$ of real type.
\end{enumerate}
\item Let $G$ be a maximal connected subgroup of $\SU n$. Then $G$ is
conjugate to one of the following:
\begin{enumerate}
\item $\SO n$;
\item $\SP k$, where $2k=n$;
\item $\mathbf{S}(\U k\times \U{n-k})$, where $1\leq k\leq n-1$;
\item $\rho(\SU p\times \SU q)$, where $pq=n$ and $p\geq 3$ and $q\geq2$,
and $\rho$ is the complex tensor product of the vector representations;
\item $\rho(G_1)$, where $G_1$ is simple and $\rho$ is a complex
irreducible representation of degree $n$ of complex type.
\end{enumerate}
\item Let $G$ be a maximal connected subgroup of $\SP n$. Then $G$ is
conjugate to one of the following:
\begin{enumerate}
\item $\U n$;
\item $\SP k\times\SP{n-k}$, where $1\leq k\leq n-1$;
\item $\rho(\SO p\times\SP q)$, where $pq=n$ and $p\geq3$ and $q\geq1$,
and $\rho$ is the real tensor product of the vector representations;
\item $\rho(G_1)$, where $G_1$ is simple and $\rho$ is a complex
irreducible representation of degree $2n$ of quaternionic type.
\end{enumerate}
\end{enumerate}
\end{thm}

Recall that a symmetric space of compact type $X=L/G$ is said to
be of \emph{inner type} if $\rk L=\rk G$; otherwise, $X$ is said
to be of \emph{outer type} (compare Theorem~8.6.7 on p.~255
in~\cite{W}). Moreover, the isotropy representation of $X$ is
absolutely irreducible if and only if $X$ is
non-Hermitian.
The following lemma implies that the isotropy representations of
symmetric spaces of semisimple type that have vanishing
homogeneity rank are precisely those coming from non-Hermitian
symmetric spaces of inner type.

\begin{lem}\label{lem:inner}
Let $(G,V)$ be the isotropy representation of
a symmetric space of compact type $X=L/G$.
Then $\homrk GV=0$ if and only if $\rk G=\rk L$.
\end{lem}

\Pf Let $\Ll=\Lg+V$ be the Cartan decomposition of $X$
with respect to the involution, where $\Ll$ and $\Lg$
respectively denote the Lie algebras of $L$ and $G$.
Let $\La\subset V$ be a maximal Abelian subspace.
By the structural theory of symmetric spaces,
it is known that the dimension of $\La$ is equal to
the cohomogeneity of $(G,V)$, and that the centralizer
$\Lm$ of $\La$ in $\Lg$ is the Lie algebra of a
principal isotropy subgroup of $(G,V)$.
It follows that $\homrk GV=\rk\Lg-\rk\Lm-\dim\La$.
Let $\Lt$ be a Cartan subalgebra of $\Lm$.
Then it is easily seen that $\Lt+\La$ is a
Cartan subalgebra of $\Ll$. Now
$\rk\Lm=\dim\Lt$, $\rk\Ll=\dim\Lt+\dim\La$, and hence
$\homrk GV=\rk\Lg-\rk\Ll=\rk G-\rk L$ which proves our thesis. \EPf

\medskip

We will also use Theorem~1.3 of~\cite{Pu} which, for convenience of the
reader, we restate here.

\begin{thm}[P\"uttmann]\label{thm:Pu}
Let $(G,M)$ be an isometric action of the compact Lie group $G$ on a
Riemannian manifold $M$. Then, for any $x\in M$, we have
\[ \dim \nu_x(G\cdot x)^{G_x} \leq \cohom GM -(\rk G_x -
\rk G_\mathrm{princ}), \]
where $\nu_x(G\cdot x)$ denotes the normal space to the orbit
$G\cdot x$ at $x$, $\nu_x(G\cdot x)^{G_x}$ denotes the fixed point subspace
of $G_x$ in $\nu_x(G\cdot x)$, and $G_\mathrm{princ}$ is a principal
isotropy subgroup of $(G,M)$.
\end{thm}

The following proposition implies that orbit-equivalent actions
have the same homogemeity rank.

\begin{prop} Let $(G,M)$ be a smooth action.
If $G'$ is a closed subgroup of $G$, and $G$ and $G'$ have the same orbits
in $M$, then $\homrk {G'}M = \homrk GM$.
\end{prop}

\Pf It is clearly enough to prove that if $G'$ and $G$ act
transitively on the same manifold $M$, then $\homrk {G'}M = \homrk GM$.
If we represent $M = G/H = G'/H'$
for suitable closed subgroups $H\subset G$ and $H'\subset G'$,
then we claim that
$$\rk G - \rk H = \rk G' - \rk H'.$$
This follows from the fact that, given a homogeneous space $M=G/H$
with $G$ compact, the number $\chi_\pi(M) := \rk H  - \rk G$
is a homotopy invariant of $M$ (see~\cite{O},
p.~207) \EPf

\medskip

Finally, we state the following
direct consequences of the definition of homogeneity rank,
which we shall repeatedly use in our arguments.

\begin{rmk}\label{rmk:estimate}
\em
Let $(G,M)$ be a smooth action of a
compact Lie group $G$ on a smooth manifold $M$. Then:
\begin{enumerate}
\item[(a)] If $\homrk GM=0$, then $\dim M\leq\dim G+\rk G$.
\item[(b)] If $G'$ is a  connected closed subgroup of $G$ having the
same homogeneity rank, then $\rk G'\geq\rk G-\rk G_\mathrm{princ}$.
Moreover if $G_\mathrm{princ}$ is finite, then $G = G'$
(indeed, $G$ and $G'$ are orbit-equivalent
both with finite principal isotropy, hence $G$ and $G'$ have the same
Lie algebra).
\end{enumerate}
\end{rmk}

\section{The classification}

In this section, we apply the strategy discussed in the previous section
to classify absolutely irreducible representations
with vanishing homogeneity rank. It is enough to consider
orthogonal representations of even degree $2n$.
According to Theorem~\ref{thm:Dynkin},
the maximal connected subgroups of $\SO{2n}$ acting
absolutely irreducibly on $V=\R^{2n}$ are:
$\rho(\SO p\times\SO q)$, where $pq=2n$ and $3\leq p\leq q$,
and $\rho$ is the real tensor product of the vector representations;
$\rho(\SP p\times\SP q)$, where $4pq=2n\neq4$,
and $\rho$ is the quaternionic tensor product of the vector representations;
and $\rho(G_1)$, where $G_1$ is simple and $\rho$ is a real form of
a complex irreducible representation of degree $2n$ of real type.

\subsection{The case of $\rho(\SO p\times\SO q)$
and its maximal subgroups}

Here $pq=2n$ and $3\leq p\leq q$.
We have that $\rho$ is the isotropy
representation of the symmetric space $\SO{p+q}/\SO p\times\SO q$,
and $\rk\SO{p+q}=\rk\SO p+\rk\SO q$ because not both of $p$,
$q$ are odd. It follows from Corollary~\ref{cor:s-repr}
that this is an example.

Next we must investigate maximal connected subgroups $G$ of
$\rho(\SO p\times\SO q)$. We shall consider separately
three cases which cover all the possibilities.

\subsubsection{$G = \rho(G_1)$, where $G_1= K \times \SO q$,
and $K \subset\SO p$ is a maximal connected
subgroup.}\label{sec:k-sop}

Set $\hat G_1=\SO p\times\SO q$, $\hat G=\rho(\hat G_1)$.
There is a $\hat G_1$-regular point
$x\in V$ whose connected principal isotropy subgroup is given by
$\hat G_{1x}=\SO{q-p}\subset\SO q$.
The isotropy subgroup of $G_1$ at $x$ is the intersection
$\hat G_{1x}\cap G_1$, and its connected component is $\SO{q-p}$.
If $G$ has vanishing homogeneity rank on $V$, then Theorem~\ref{thm:Pu}
applied to $x$ gives
\begin{eqnarray*}
\dim\nu_x(Gx)^{G_x} &\leq& \cohom{G}{V} - (\rk G_x - \rk G_\mathrm{princ} ) \\
                    &=& \rk G - \rk G_x \\
                    &=&\rk K + \rk\SO q -\rk\SO{q-p} \\
                    &=&\rk K + \rk\SO p,
\end{eqnarray*}
where we have used that not both of $p$, $q$ are odd. Note that
$\dim \nu_x(\hat Gx) = p$ and $\dim \nu_x(Gx) = p +\dim\SO p -\dim K$.
It is clear that
\[ \nu_x(Gx)^{G_x} \supset \nu_x(\hat Gx), \]
since $\nu_x(Gx) \supset \nu_x(\hat Gx)$,
$G_x\subset \hat G_x$ and $x$ is $\hat G$-regular.
It follows that
\[ \dim\nu_x(Gx)^{G_x} \geq p. \]
Combining with the above we get that
\[ \rk K + \rk\SO p \geq p \geq 2\rk\SO p, \]
and therefore
\[ \rk K = \rk\SO p. \]
By the classification of maximal subgroups of maximal rank of $\SO p$,
(see, for example, section~8.10 in~\cite{W}),
since $K$ is irreducible and of real type on $\R^p$, we
must have $K=\SO p$.

\subsubsection{$G = \rho(G_1)$, where $G_1=\SO p\times K$,
and $K\subset\SO q$ is a maximal connected subgroup.}

According to Theorem~\ref{thm:Dynkin}, we need
to consider three cases.

\paragraph{(a) $K=\mu(K_1)$, where $K_1$ is simple
and $\mu$ is an absolutely irreducible representation
of degree $q$.}
Of course we need only to consider representations $\mu$
such that $\mu(K_1)$ is a proper subgroup of $\SO q$.

We may assume $p<q$. Remark~\ref{rmk:estimate}(a) gives that
$p^2 < pq \leq \frac{p(p-1)}2 + \left[\frac p2\right] + r$,
where $r = \dim K_1 + \rk K_1$, and $[x]$ denotes the greatest
integer contained not exceeding $x$. This implies that
\begin{equation}\label{2}
 p^2 < 2r,
\end{equation}
and that
\begin{equation}\label{3}
p^2-2qp+2r \geq 0.
\end{equation}
Equation~(\ref{3}) and $p\geq3$ then imply that
\begin{equation}\label{4}
q \leq \frac r3 + \frac32.
\end{equation}
Let $s$ be the minimal degree of an absolutely irreducible representation
of $K$ and such that its image is not the full
$\SO s$.
Then
\begin{equation}\label{5}
s > \frac r3 + \frac32
\end{equation}
is a sufficient condition for $(G,V)$ not to be of vanishing
homogeneity rank. We next run through the possibilities for $K_1$.
\begin{itemize}
\item $K_1=\SU m$, where $m\geq2$. Here $r=m^2+m-2$. If $m=2$,
then $s=5$, and~(\ref{5}) holds. If $m\geq3$, then $s=m^2-1$
(realized by the adjoint representation), and~(\ref{5}) holds.
\item $K_1=\SP m$, where $m\geq2$. Here $r=2m^2+2m$. If $m=2$,
then $s=10$ (realized by the adjoint representation) and~(\ref{4})
holds. If $m\geq3$ then $s=2m^2-m-1$ (realized by the second
fundamental representation) and~(\ref{4}) holds. \item $K_1=\Spin
m$, where $m\geq 7$. Here $r=\frac{m(m-1)}2+\left[\frac
m2\right]$. All irreducible representations of real type
violate~(\ref{4}), except possibly the (half-)spin representations
of $\B{4k-1}$, $\B{4k}$, $\D{4k}$. These have respectively
$q=2^{4k-1}$, $2^{4k}$, $2^{4k-1}$. Condition~(\ref{4}) is
respectively
\[ 3\cdot2^{4k} \leq 64k^2-16k+9, \]
\[ 3\cdot2^{4k+1} \leq 64k^2+16k+9, \]
\[ 3\cdot2^{4k} \leq 64k^2+9. \]
The only cases that survive are $\B3$ and $\D4$.
In the case of $\B3$ we have $q=8$ and $r=24$.
Then~(\ref{2}) implies that $p=3$, $4$, $5$, $6$.
Next we use~(\ref{3}) to get rid of $p=5$, $6$.
We end up with $p=3$ and $p=4$, and this gives the
admissible cases $(\SO3\times\Spin7,\R^3\otimes\R^8)$ and
$(\SO4\times\Spin7,\R^4\otimes\R^8)$, but note that the
first one of these is orbit-equivalent to
$(\SO3\times\SO8,\R^3\otimes\R^8)$.
In the case of $\D4$, we have that $\mu(\Spin8)=\SO8$,
 and we rule this out.
\item $K_1$ is an exceptional group. Here~(\ref{5}) holds
in each case, so there are no examples, see the table below.
\[\begin{array}{|c|c|c|}
\hline
 K_1 & r & s \\
\hline
\G & 16 & 7 \\
\F & 56 & 26 \\
\E6 & 84 & 78 \\
\E7 & 140 & 133 \\
\E8 & 256 & 248 \\
\hline
\end{array}\]
\end{itemize}

\paragraph{(b) $K=\mu(\SO k\times\SO l)$,
where $3\leq k\leq l$ and $q=kl$,
and $\mu$ is the real tensor product of the vector representations.}
Here $r = \frac{k^2+l^2}2 + \underbrace{\left[\frac k2\right]-\frac k2
+ \left[\frac l2\right]-\frac l2}_{=\theta}$.
Note that $-1\leq\theta\leq0$. Then~(\ref{4}) is
$k^2+l^2-6kl+2\theta+9 \geq0$.
Set $m=l-k\geq0$. Then $m^2-4km-4k^2+9+2\theta \geq 0$.
This implies that
\begin{equation}\label{6}
 2k+ \sqrt{8k^2-9-2\theta}\leq m.
\end{equation}
If $pk\leq l$, since the action of
$\SO p\times\mu(\SO k\times\SO l)\subset\SO p\times\SO q$
on $\R^p\otimes\R^q$ is the same thing as the action of
$\mu'(\SO p\times\SO k)\times\SO l\subset\SO{pk}\times\SO l$
on $\R^{pk}\otimes\R^l$, where $\mu'$ is the real tensor
product of the vector representations,
this case has already been considered
in section~\ref{sec:k-sop}. So now we assume that
\begin{equation}\label{7}
pk>l.
\end{equation}
Note that $q>\sqrt{2r}$.
This implies via~(\ref{3}) that $3\leq p\leq q-\sqrt{q^2-2r}$.
Combining this with~(\ref{7}) we have $l<kq-k\sqrt{q^2-2r}$
and then $l^2(1-2k^2)+2k^2r>0$.
Substituting the value of $r$ we get
$k\leq l < k \sqrt{\frac{k^2+2\theta}{k^2-1}}$.
We deduce that $\theta=0$ and
\begin{equation}\label{8}
 0\leq m < k\left(\frac{k}{\sqrt{k^2-1}}-1\right).
\end{equation}
Now~(\ref{6}) and~(\ref{8}) combined imply that
$3k+\sqrt{8k^2-9} < \frac{k^2}{\sqrt{k^2-1}}$,
which is impossible for $k\geq3$.

\paragraph{(c) $K=\mu(\SP k\times\SP l)$, where $q=4kl\neq4$,
and $\mu$ is the quaternionic tensor product of the vector representations.}
We postpone this case to section~\ref{sec:so-sp-sp}.

\subsubsection{$G=\{(x,\sigma(x)):x\in\SO p \}$, where
$p=q$ and $\sigma$ is an automorphism of $\SO p$.}

Here Remark~\ref{rmk:estimate}(a)
immediately implies that $(G,V)$ cannot have vanishing
homogeneity rank.

\subsection{The case of $\rho(\SP p\times\SP q)$ and its
maximal subgroups}

Here $4pq=2n\neq4$ and $p\leq q$.
We have that $\rho$ is the isotropy
representation of the symmetric space $\SP{p+q}/\SP p\times\SP q$,
and $\rk\SP{p+q}=\rk\SP p+\rk\SP q$.
It follows from Corollary~\ref{cor:s-repr}
that this is an example.

Next we must investigate maximal connected subgroups $G$ of
$\rho(\SP p\times\SP q)$. We shall consider three cases
separately which cover all the possibilities.

\subsubsection{$G = \rho(G_1)$, where
$G_1 = K \times \SP q$, and $K \subset\SP p$
is a maximal connected subgroup.}\label{sec:k-x-spq}

Set $\hat G_1=\SP p\times\SP q$. There is a $\hat G_1$-regular point
$x\in V$ whose principal isotropy subgroup is given by $\hat
G_{1x}=\SP1^p\times\SP{q-p}$. Let $\mathcal O$ be the orbit of $\hat
G_1$ through $x$. A point $y\in\mathcal O$ is $G_1$-regular in
$\mathcal O$ if and only if it is $G_1$-regular in $V$; moreover the
isotropy subgroup $G_{1y}$ is given by the intersection of $G_1$ with a
suitable conjugate of $\hat G_{1x}$ in $\hat G_1$. This means that a
principal isotropy subgroup of $G_1$ contains a subgroup isomorphic
to $\SP{q-p}$. If $(G,V)$ has vanishing homogeneity rank, it
then follows that
\[ \dim V=4pq\leq\dim K+\rk K+2q^2+2q-(2(q-p)^2+2(q-p)), \]
hence
\begin{equation}\label{11}
 2p^2-2p\leq \dim K+\rk K.
\end{equation}
According to Theorem~\ref{thm:Dynkin}, there are two cases
to be considered. But if $K$ is of the form $\mu(\SO k\times\SP l)$,
where $\mu$ is the real tensor product of the vector representations,
we refer to section~\ref{sec:so-sp-sp}. So we can assume
that $K$ is of the form $\mu(K_1)$, where $K_1$ is simple and
$\mu$ is a complex irreducible representation of degree $2p$
of quaternionic type.
Let $s$ be half the minimal
degree of a complex irreducible representation
of quaternionic type of $K_1$ and such that its image is not the full
$\SP s$. The list of the values of $s$ for each compact
simple group $K_1$ is given by the following table
(groups not appearing in the table do not admit quaternionic
representations):

\setlength{\extrarowheight}{0.15cm}

\[\begin{array}{|c|c|}
\hline
 K_1 & s \\
\hline\hline
\SU{4a+2},\ a\geq1& \frac12\binom{4a+2}{2a+1} \\ \hline
\Spin{8a+3},\ a\geq1& 2^{4a} \\ \hline
\Spin{8a+4},\ a\geq1& 2^{4a} \\ \hline
\Spin{8a+5},\ a\geq1& 2^{4a+1} \\ \hline
\SP a,\ a\geq3 & \frac a3(2a^2-3a-2) \\ \hline
\SP 2 & 8 \\ \hline
\SP 1 & 2 \\ \hline
\E7 & 28 \\
\hline
\end{array}\]
It is now easy to see that~(\ref{11}) implies
$K_1=\SP1$ and $p=2$, and then the only admissible
case is $G_1=\SP1\times\SP q$, where $\SP1\subset\SP2$ via
the irreducible representation of degree~$4$.

\subsubsection{$G = \rho(G_1)$, where $G_1=\SP p\times K$,
and $K\subset\SP q$ is a maximal connected subgroup.}\label{sec:spp-x-k}

According to Theorem~\ref{thm:Dynkin}, we need to consider
two cases.

\paragraph{(a) $K=\mu(K_1)$, where $K_1$ is simple
and $\mu$ is a complex irreducible representation
of degree $2q$ of quaternionic type.}
We may assume $p<q$. Remark~\ref{rmk:estimate}(a) gives that
$4p^2<4pq\leq 2p^2+2p+r$,
where $r = \dim K_1 + \rk K_1$ (note that $r\geq4$). This implies that
\[ p < \frac{1+\sqrt{1+2r}}2\quad\mbox{and}\quad
2p^2+2p(1-2q)+r\geq0. \]
From this we get that
\begin{equation}\label{13}
 q \leq \frac r4 + 1,
\end{equation}
and
\begin{equation}\label{14}
\mbox{if $q\geq \frac{1+\sqrt{1+2r}}2$,
then $p\leq q-\frac12-\frac12\sqrt{(2q-1)^2-2r}$}.
\end{equation}
Running through the compact simple groups
$K_1$ that admit quaternionic representations
(see table in section~\ref{sec:k-x-spq})
and using~(\ref{13}) and~(\ref{14}), we get the following admissible
cases:
$K_1=\SP1$, $p=1$, $q=2$;
$K_1=\SP3$, $p=1$, $q=7$;
$K_1=\Spin{11}$, $p=1$, $q=16$;
$K_1=\Spin{12}$, $p=1$, $q=16$;
$K_1=\SU6$, $p=1$, $q=10$;
$K_1=\E7$, $p=1$, $q=28$.
All cases but that of $K_1=\Spin{11}$ come from isotropy representations
of symmetric spaces.

\paragraph{(b) $K=\mu(\SO k\times\SP l)$, where $q=kl$,
and $\mu$ is the real tensor product of the vector representations.}
We postpone this case to section~\ref{sec:so-sp-sp}.

\subsubsection{$G=\{(x,\sigma(x)):x\in\SP p \}$, where
$p=q$ and $\sigma$ is an automorphism of $\SP p$.}

Here Remark~\ref{rmk:estimate}(a) immediately implies that $(G,V)$
can have vanishing homogeneity rank only if $p=1$, so this case is
out.

\subsection{The case of $\rho(G_1)$}

Here $G_1$ is a compact simple Lie group and $\rho$ is an
absolutely irreducible representation of $G_1$ of degree
$2n$. Remark~\ref{rmk:estimate}(a) says that
$2n\leq\dim G_1 + \rk G_1$. In particular, this implies
that $2\dim G_1\geq 2n - 2$, so we can use Lemma~2.6 in~\cite{K}
to deduce that $(G,V)$ is orbit equivalent to the isotropy representation
of a symmetric space.

\subsection{The case of $\rho(\SO m\times\SP p\times\SP q)$,
where $\rho$ is the real and quaternionic
tensor products of the vector representations}
\label{sec:so-sp-sp}

Here $2n=4mpq$, $m\geq3$ and $p\leq q$.
By direct computation or using Theorem~1.1 in~\cite{HH},
we see that:
\begin{enumerate}
\item[(i)] if $m\geq 4pq+2$, then the connected principal
isotropy is given by $\SO{m-4pq}$;
\item[(ii)] if $q\geq mp+1$, then the connected principal
isotropy is given by $\SP{q-mp}$;
\item[(iii)] in all other cases the connected principal
isotropy is trivial.
\end{enumerate}
In case (i) the condition of vanishing homogeneity rank reads
\[ 4p^2q^2=p^2+q^2+p+q\leq 2p^2+2q^2, \]
and this implies $p=q=1$.
In case (ii) we have
\[ 4p^2(1-m^2)+4p(1+m)+m^2-m+2\left[\frac m2\right]=0. \]
If $m=2l$, then we have $p^2(1-4l^2)+p(1+2l)+l^2=0$,
which implies that $1+2l$ divides $l^2$, impossible.
If $m=2l+1$, then we have $(4p^2-1)l=2p$,
which is impossible.
In case (iii), we have the equation
\[ 8mpq=m^2-m+2\left[\frac m2\right]+4p^2+4p+4q^2+4q. \]
If $m=2l$, this reads
\begin{equation}\label{9}
 l^2 - 4pql +p^2+q^2+p+q=0,
\end{equation}
subject to the constraints
\[ \frac q{2p}\leq l \leq 2pq,\quad p\leq q,\quad l\geq2, \]
while if $m=2l+1$, we have
\begin{equation}\label{10}
 l^2 -(4pq-1)l+ p^2+q^2-2pq+p+q=0,
\end{equation}
subject to the constraints
\[ \frac q{2p}-\frac12 \leq l \leq 2pq,\quad p\leq q,\quad l\geq1. \]

Consider first equation~(\ref{9}). It can be solved in~$l$ to yield
$l=2pq\pm\sqrt\Delta$,
where
$\Delta=4p^2q^2-p^2-q^2-p-q$. If $l=2pq+\sqrt\Delta$, using the fact
that $l\leq 2pq$ we have $\Delta=0$ and then $l=2pq\leq p+q$,
which gives $p=q=1$, and then $l=2$, $m=4$.
If $l=2pq-\sqrt\Delta$, then $\frac q{2p}\leq l$
implies that
\[ q^2(4p^2-1)-4p^2q-4p^2(p^2+p)\leq 0, \]
and therefore
\[ p\leq q \leq \frac{2p^2+2p\sqrt{p(4p^3+4p^2-1)}}{4p^2-1}<p+2. \]
It then follows that we only need to consider
the possibilities $q=p$ and $q=p+1$.
If $q=p$, then we have
\[ 2\leq l=2p^2-\sqrt{4p^4-2p^2-2p}\leq2, \]
so that $l=2$, $p=q=1$.
If $q=p+1$, then
\[ 2\leq l=2p^2+2p-\sqrt{4p^4+8p^3+2p^2-4p-2}<2, \]
which is impossible.

Next we consider equation~(\ref{10}).
Here it is useful to note that Remark~\ref{rmk:estimate}(b) applied
to $G=\SO m\times\SP p\times\SP q$ viewed as a subgroup
of $\SO m\times\SO{4pq}$ gives the
extra condition
\begin{equation}\label{12}
\left[\frac m2\right] \leq p+q.
\end{equation}
Equation~(\ref{10}) can be solved in $l$ to yield
$l=\frac{4pq-1\pm\sqrt{\Delta_1}}2$,
where $\Delta_1=16p^2q^2-4p^2-4p-4q^2-4q+1$.
If $l=\frac{4pq-1+\sqrt{\Delta_1}}2$, then~(\ref{12})
implies that
\[ 2pq\leq l \leq p+q \leq 2q, \]
which gives $p=q=1$, $\Delta_1=1$, $l=2$ and $m=5$.
If $l=\frac{4pq-1-\sqrt{\Delta_1}}2$, then the inequality
$\frac q{2p}-\frac12 \leq l$ implies that
\[ q^2(4p^2-1)-4p^2q-4p^4-4p^3+p^2\leq0, \]
which in turn implies
\[ p\leq q\leq \frac{2p^2+p\sqrt{16p^4+16p^3-4p^2-4p+1}}{4p^2-1}<p+2. \]
Therefore we need only to consider the cases $q=p$ and $q=p+1$.
If $q=p$, then
\[ l=\frac{4p^2-1-\sqrt{16p^4-8p^2-8p+1}}2 < 2. \]
This gives $l=1$, $m=3$, $p=q=1$. If $q=p+1$, then
\[ l=2p^2+2p-\frac12-\frac12\sqrt{16p^4+32p^3+8p^2-16p-7}<1, \]
which is impossible.

\subsection{The examples and their subgroups}
\label{sec:ex-sbg}

In this section we show that all candidates $G\subset\SO{2n}$ found in the
previous sections are actually examples of groups acting
absolutely irreducibly on $V=\R^{2n}$ with vanishing homogeneity rank, and
these groups do not admit subgroups with the same property.
This will complete the proof of Theorem~\ref{thm:main}.

We first examine the representations that are not orbit-equivalent to
isotropy representations of non-Hermitian symmetric spaces of inner type.
We have three candidates:
\begin{enumerate}
\item $G = \SP1\times\SP q$ ($q\geq 2$) acting on
$V = S^3({\mathbf C}^2)\otimes_{\mathbf H}{\mathbf C}^{2q}\cong\R^{8q}$;
\item $G = \SO4\times\Spin7$ acting on $V = \R^4\otimes\R^8\cong\R^{32}$, where $\Spin7$ acts on $\R^8$ via the spin representation;
\item $G = \SP1\times \Spin{11}$ acting on $V = {\mathbf C}^2\otimes_{{\mathbf H}}{\mathbf C}^{32}\cong \R^{64}$, where
$\Spin{11}$ acts on ${\mathbf C}^{32}$ via the spin representation.
\end{enumerate}
We now show that in each case the representation has vanishing homogeneity
rank. Indeed, in case 1 we have that a
connected principal isotropy is given by $\SP{q-2}$
(see~\cite{GT}, Proposition~7.12),
therefore the cohomogeneity is three and
the homogeneity rank vanishes.
In case 2, a connected isotropy subgroup is trivial.
This can be seen by selecting a pure tensor $v\otimes w$
with $v\in\R^4$ and $w\in\R^8$ and computing the connected isotropy,
which is $\SO3\times\G$; then the slice
representation is given by $\R\oplus \R^3\otimes \R^7$;
starting again with this new representation, we eventually
come up with a trivial isotropy. Therefore the cohomogeneity is five and
the homogeneity rank vanishes.
In case 3, we also have trivial connected principal isotropy and vanishing
homogeneity rank. Indeed, if $v\in{\mathbf C}^{32}$ is a highest weight
vector for the spin representation of $\Spin{11}$,
then the subgroup $H\subset \Spin{11}$ defined by
$H = \{g\in \Spin{11}:\;g\cdot v \in {\mathbf C}^*\cdot v\}$ is
given by $\U5$.
Now if $p:\SP1\times\Spin{11}\to \Spin{11}$ is the projection, then
$$ p((\SP1\times\Spin{11})_v) = \{g\in\Spin{11}:\;g\cdot v \in \SP1\cdot v\} \supset H.$$
Since $H$ is maximal in $\Spin{11}$, we get that $(\SP1\times\Spin{11})_v$
is given by ${\mathbf T}^1\cdot \SU5$,
where ${\mathbf T}^1$ sits diagonally in the product of a suitable
maximal torus in $\SP1$ and the center of $H$.
From this we see that the slice representation at $v$ is given
by $\R \oplus {\mathbf C}^5 \oplus \Lambda^2{\mathbf C}^5$
and the connected principal isotropy is trivial. The cohomogeneity is six
and the homogeneity rank vanishes.

We now examine subgroups of the previous examples.
In case~1, a maximal subgroup of $G$ leaving no complex structure
on $V$ invariant is of the form
$G' = \SP1\times K$, where $K\subset \SP{q}$ is maximal.
Since $\SP2\times K$ does not have vanishing homogeneity rank on $V$
by the results of section~\ref{sec:spp-x-k}, and $G'\subset\SP2\times K$,
we have that $G'$ does not have vanishing homogeneity rank on $V$.
In cases~2 and~3, $G$ admits no proper subgroups acting with vanishing
homogeneity rank because the connected principal isotropy is trivial and then
we may apply Remark~\ref{rmk:estimate}(b).

We finally consider the representations $(G,V)$ that are orbit-equivalent to
isotropy representations of non-Hermitian symmetric spaces of inner type,
and we classify the subgroups $G'\subset G$
which still act absolutely irreducibly on $V$ with
vanishing homogeneity rank.
In the following table we list the representations $\rho$ which need to be
examined; we denote by $c$ the cohomogeneity of $\rho$,
by $d$ the dimension of $V$, and
by $[[W]]$ a real form of the $G$-module $W$.

\setlength{\extrarowheight}{0.15cm}

\[\begin{array}{|c|c|c|c|c|c|}

\hline \textsl{Case} &G & \rho & c & d & \dim G_{\text{princ}} \\
\hline 1&  \SP1\cdot\SU6 & \CC^2\otimes_{\mathbf{H}}\Lambda^3\CC^6  & 4 & 40 & 2\\
\hline 2 & \SP1\cdot\Spin{12} & \CC^2\otimes_{\mathbf{H}}{\text{(half-spin)}} &4 & 64 & 9 \\
\hline 3 & \SP1\cdot\E7 & \CC^2\otimes_{\mathbf{H}}\CC^{56}  &4 & 112 & 28\\
\hline 4&  \SP1\cdot\SP3 & \CC^2\otimes_{\mathbf{H}}\CC^{14}&4 & 28 & 0\\
\hline 5&  \SO4 & S^3(\CC^2)\otimes_{\mathbf{H}}\CC^2 &2 & 8 & 0\\
\hline 6 & \Spin{16} & {\text{half-spin}}   &8& 128 & 0\\
\hline 7& \SU8  &[[\Lambda^4\CC^8]] & 7 & 70 & 0\\
\hline 8 & \SO3\times\Spin7 & \R^3\otimes\R^8   &3& 24 & 3\\
\hline 9&  \Spin7 & {\text{spin}} &1& 8 & 14\\
\hline 10&  \Spin9 & {\text{spin}} &1& 16 & 21\\

 \hline
\end{array}\]

Cases~4 through~7 can be dealt with using Remark~\ref{rmk:estimate}(b).
We consider case~1. If $G'$ is a maximal subgroup of $G$,
we may assume that $G'$ is of the form $G' = \SP1\cdot G''$,
where $G''$ is maximal in $\SU6$, since $G'$ does not
leave any complex structure invariant. Now Remark~\ref{rmk:estimate}(a)
implies that $\dim G'' + \rk G'' \geq 36$, and $\rk G''\leq 5$ implies
$\dim G'' \geq 31$, so that
$\dim \SU6/G'' \leq 4$. If $G''$ is a proper subgroup of
$\SU6$, then the left action of $\SU6$ on
$\SU6/G''$ is almost effective because $\SU6$ is simple.
Therefore $\dim\SU6$ is less than the dimension
of the isometry group of $\SU6/G''$, which is at most $10$,
but this is a contradiction.
Hence $G''=\SU6$.

In case~2, again we can assume that $G'$ is of the form
$G' = \SP1\cdot G''$, where $G''$ is maximal in $\Spin{12}$. We have
$\dim G'' + \rk G'' \geq 60$; since $G'$ is supposed to act
absolutely irreducibly on $V$, its rank is not maximal by a
theorem of Dynkin
(see Theorem~7.1, p.~158 in~\cite{D2}),
and therefore $\dim G''\geq 55 = \dim \Spin{11}$.
It is known that a subgroup of $\Spin{n}$ of dimension greater or equal
to $\dim \Spin{n-1}$ is conjugate to the standard
$\Spin{n-1}\subset \Spin{n}$ if $n\neq 4$, $8$
(see e.g.~\cite{Ko}, p.~49).
So, $G' = \SP1\cdot\Spin{11}$,
which is indeed an example with trivial connected principal isotropy
by the discussion above.

In case~3, using the same argument as in case~2, we see that
$G'=\SP1\cdot G''$, where $\dim G''\geq 102$. An inspection of the list of
all maximal subalgebras of $\E7$
(see Table~12, p.~150 and Theorem~14.1, p.~231
in~\cite{D2}) shows that there is no such proper subgroup.

In case~8, a maximal subgroup $G'$ acting absolutely
irreducibly on $V$ must be of the form $G' = \SO3\times K$, where
$K\subset\Spin7$ is maximal; arguing as above,
we see that $\dim K \geq 18$, so that
$\dim (\Spin7/K)\leq 3$ and this is impossible, because $\Spin7$ is simple.

In case~9, let $K\subset \Spin7$ be a maximal subgroup acting
absolutely irreducibly on $\R^8$. Since $K$ cannot have maximal rank
as above, and using Theorem~\ref{thm:Dynkin},
we see that $K$ must be simple of rank at most two and
it must admit an irreducible representation of degree~$7$ and of real type.
Moreover, by Remark~\ref{rmk:estimate}(a), we have $\dim K \geq 6$,
hence $\rk K = 2$, and a direct inspection
of all such simple groups shows that none of them but $\G$ admits an
irreducible representation of degree~$7$. But
$\G$ does not admit an irreducible representation of degree~$8$.

In case~10, we consider a maximal subgroup $K$ of $\Spin9$ acting
absolutely irreducibly on $\R^{16}$.
This means that $\rk K \leq 3$ and $\dim K\geq 13$. Looking at the list
of all maximal subgroups of $\Spin9$, we see that we can suppose $K$ to
be simple and to act irreducibly on $\R^9$, via
the embedding $K\subset \Spin9$. Therefore $K$ must be one of
$\G$, $\SU4$, $\Spin7$ or $\SP3$,
but we immediately see
that none of these groups admits an irreducible representation of
degree $9$ and of real type.

This finishes the proof of Theorem~\ref{thm:main}.


\providecommand{\bysame}{\leavevmode\hbox to3em{\hrulefill}\thinspace}

\end{document}